\documentclass[runningheads]{llncs}
\usepackage[breaklinks, colorlinks, linkcolor=blue, citecolor=blue, urlcolor=black]{hyperref}
\usepackage{color}

\usepackage{graphicx}
\usepackage{tabularx}
\usepackage{tabularcalc}
\usepackage{booktabs}
\usepackage{multirow}
\usepackage{algorithm}
\usepackage{algpseudocode}
\usepackage{tikz}
\usepackage{forest}
\usepackage{pgfplots}
\usepackage{xltabular}
\usetikzlibrary{arrows.meta}
\usepackage{amsmath}
\usepackage{amssymb}
\usepackage{xspace}
\usepackage{xcolor}
\usepackage{varioref}
\usepackage{cleveref}
\usepackage{aligned-overset}
\usepackage[shortlabels,inline]{enumitem}
\usepackage{boxedminipage2e}
\usepackage{makecell}

\definecolor{pbblue}{RGB}{68,119,170}
\definecolor{pbcyan}{RGB}{102,204,238}
\definecolor{pbgreen}{RGB}{34,136,51}
\definecolor{pbred}{RGB}{238,102,119}

\usepackage{ifthen}

\usepackage{ifthen}

\provideboolean{detectedSTOC}
\provideboolean{detectedFOCS}
\provideboolean{detectedElsevier}
\provideboolean{detectedNOW}
\provideboolean{detectedLMCS}
\provideboolean{detectedIEEE}
\provideboolean{detectedPoster}
\provideboolean{detectedSIAM}
\provideboolean{detectedLNCS}
\provideboolean{detectedACM}
\provideboolean{detectedACMconf}
\provideboolean{detectedSigplanconf}
\provideboolean{detectedToC}
\provideboolean{detectedLIPIcs}
\provideboolean{detectedAAAI}
\provideboolean{detectedIJCAI}
\provideboolean{detectedCompCplx}
\provideboolean{detectedEasyChair}
\provideboolean{detectedArticle}
\provideboolean{detectedReport}
\provideboolean{detectedThesis}

\makeatletter

\@ifclassloaded{sig-alternate}
{\setboolean{detectedSTOC}{true}}
{\setboolean{detectedSTOC}{false}}

\@ifclassloaded{elsarticle}
{\setboolean{detectedElsevier}{true}}
{\setboolean{detectedElsevier}{false}}

\@ifclassloaded{now}
{\setboolean{detectedNOW}{true}}
{\setboolean{detectedNOW}{false}}

\@ifclassloaded{lmcs}
{\setboolean{detectedLMCS}{true}}
{\setboolean{detectedLMCS}{false}}

\@ifclassloaded{IEEEtran} {
  \setboolean{detectedIEEE}{true}
  \ifCLASSOPTIONconference {
    \setboolean{detectedFOCS}{true}
  }
  \else {
    \setboolean{detectedFOCS}{false}
  }
  \fi
}
{
  \setboolean{detectedFOCS}{false}
  \setboolean{detectedIEEE}{false}
}

\@ifclassloaded{siamart171218}
{\setboolean{detectedSIAM}{true}}
{\setboolean{detectedSIAM}{false}}

\@ifclassloaded{llncs}
{\setboolean{detectedLNCS}{true}}
{\setboolean{detectedLNCS}{false}}

\@ifclassloaded{acmsmall}
{\setboolean{detectedACM}{true}}
{\setboolean{detectedACM}{false}}

\@ifclassloaded{acmart}
{\setboolean{detectedACMconf}{true}}
{\setboolean{detectedACMconf}{false}}

\@ifclassloaded{sigplanconf}
{\setboolean{detectedSigplanconf}{true}}
{\setboolean{detectedSigplanconf}{false}}

\@ifclassloaded{toc}
{\setboolean{detectedToC}{true}}
{\setboolean{detectedToC}{false}}

\@ifclassloaded{lipics}
{\setboolean{detectedLIPIcs}{true}}
{\@ifclassloaded{lipics-v2019}
  {\setboolean{detectedLIPIcs}{true}}
  {\@ifclassloaded{oasics-v2019}
    {\setboolean{detectedLIPIcs}{true}}
    {\setboolean{detectedLIPIcs}{false}}
  }
}

\@ifclassloaded{cc}
{\setboolean{detectedCompCplx}{true}}
{\setboolean{detectedCompCplx}{false}}

\@ifclassloaded{easychair}
{\setboolean{detectedEasyChair}{true}}
{\setboolean{detectedEasyChair}{false}}

\@ifpackageloaded{aaai}
{\setboolean{detectedAAAI}{true}}        
{\@ifpackageloaded{aaai18}
  {\setboolean{detectedAAAI}{true}}
  {\@ifpackageloaded{aaai20}       
    {\setboolean{detectedAAAI}{true}}
    {\setboolean{detectedAAAI}{false}}}}

\@ifpackageloaded{ijcai18}
{\setboolean{detectedIJCAI}{true}}        
{\@ifpackageloaded{ijcai19}
  {\setboolean{detectedIJCAI}{true}}        
  {\setboolean{detectedIJCAI}{false}}}

\@ifclassloaded{sciposter}
{\setboolean{detectedPoster}{true}}
{\setboolean{detectedPoster}{false}}

\@ifclassloaded{article}
{\setboolean{detectedArticle}{true}}
{\setboolean{detectedArticle}{false}}

\makeatother

\ifthenelse{\not \isundefined{\examen} 
  \and \not \isundefined{\disputationsdatum} 
  \and \not \isundefined{\disputationslokal}}   
  {\setboolean{detectedThesis}{true}}
  {\setboolean{detectedThesis}{false}}

\ifthenelse{\boolean{detectedArticle} \or \boolean{detectedThesis}
  \or \boolean{detectedSTOC} \or \boolean{detectedFOCS}
  \or \boolean{detectedSIAM} \or \boolean{detectedIEEE}
  \or \boolean{detectedPoster}}
{\setboolean{detectedReport}{false}}
{\setboolean{detectedReport}{true}}

\usepackage{xspace}
\ifthenelse
{\boolean{detectedElsevier} \or \boolean{detectedSIAM} 
  \or \boolean{detectedLIPIcs}}
{}
{\usepackage{varioref}}

\makeatletter
\AtBeginDocument{%
  \def\doi#1{\url{https://doi.org/#1}}}
\makeatother

\DeclareMathAlphabet{\mathsfsl}{OT1}{cmss}{m}{sl}

\DeclareRobustCommand{\BibTeX}{%
  {\normalfont B\kern-.05em{\scshape i\kern-.025em b}\kern-.08em \TeX}%
}

\ifthenelse{\boolean{detectedToC}}{}
{
  \newcommand{\Q}         {\mathbb{Q}}

}

\newcommand{\ceiling}[1]{\lceil #1 \rceil}

\newcommand{\CEILING}[1]{\left \lceil #1 \right \rceil}

\newcommand{\MAXOFEXPR}[2][]{\max_{#1} \left\{ #2 \right\}}
\newcommand{\MINOFEXPR}[2][]{\min_{#1} \left\{ #2 \right\}}
\newcommand{\Maxofexpr}[2][]{\max_{#1} \bigl\{ #2 \bigr\}}
\newcommand{\Minofexpr}[2][]{\min_{#1} \bigl\{ #2 \bigr\}}

\newcommand{\MAXOFSET}[3][:]%
     {\ifthenelse{\equal{#1}{;}}%
     {\MAXOFEXPR{ #2 \,;\, #3 }}
     {\ifthenelse{\equal{#1}{:}}%
     {\MAXOFEXPR{ #2 \,:\, #3 }}
     {\max \twincommandJN{\left\{}{#2}{\left#1}{\right}{\,#3}{\right\}}}}}
\newcommand{\MINOFSET}[3][:]%
     {\ifthenelse{\equal{#1}{;}}%
     {\MINOFEXPR{ #2 \,;\, #3 }}
     {\ifthenelse{\equal{#1}{:}}%
     {\MINOFEXPR{ #2 \,:\, #3 }}
     {\min \twincommandJN{\left\{}{#2}{\left#1}{\right}{\,#3}{\right\}}}}}

\newcommand{\Maxofset}[3][:]%
     {\ifthenelse{\equal{#1}{;}}%
     {\Maxofexpr{ #2 \,;\, #3 }}
     {\ifthenelse{\equal{#1}{:}}%
     {\Maxofexpr{ #2 \,:\, #3 }}
     {\max \twincommandJN{\bigl\{}{#2}{\bigl#1}{\bigr}{\,#3}{\bigr\}}}}}
\newcommand{\Minofset}[3][:]%
     {\ifthenelse{\equal{#1}{;}}%
     {\Minofexpr{ #2 \,;\, #3 }}
     {\ifthenelse{\equal{#1}{:}}%
     {\Minofexpr{ #2 \,:\, #3 }}
     {\min \twincommandJN{\bigl\{}{#2}{\bigl#1}{\bigr}{\,#3}{\bigr\}}}}}

\DeclareMathOperator{\Expop}{E}

\ifthenelse{\boolean{detectedLMCS}}
{}
{}

\newcommand{\twincommandJN}[6]%
    {#1#2#3\vphantom{#2#5}\mspace{-2.05mu}#4.#5#6}

\newcommand{\CondExp}[2]%
    {\Expop\twincommandJN{\bigl[}{#1}{\bigl|}{\bigr}{\,#2}{\bigr]}}
\newcommand{\CONDEXP}[2]%
     {\Expop\twincommandJN{\left[}{#1}{\left|}{\right}{\,#2}{\right]}}

\newcommand{\Condprob}[3][]%
    {\Pr_{#1}\twincommandJN{\bigl[}{#2}{\bigl|}{\bigr}{\,#3}{\bigr]}}
\newcommand{\CONDPROB}[3][]%
    {\Pr_{#1}\twincommandJN{\left[}{#2}{\left|}{\right}{\,#3}{\right]}}

\newcommand{\set}[1]{\{ #1 \}}
\newcommand{\Set}[1]{\bigl\{ #1 \bigr\}}

\newcommand{\Setdescr}[3][|]%
     {\ifthenelse{\equal{#1}{;}}%
     {\Set{ #2 \,;\, #3 }}
     {\ifthenelse{\equal{#1}{:}}%
     {\Set{ #2 \,:\, #3 }}
     {\twincommandJN{\bigl\{}{#2\,}{\bigl#1}{\bigr}{\,#3}{\bigr\}}}}}
\newcommand{\SETDESCR}[3][|]%
     {\twincommandJN{\left\{}{#2\,}{\left#1}{\right}{\,#3}{\right\}}}

\newcommand{\Setdescrbrackets}[3][|]%
     {\twincommandJN{\bigl[}{#2}{\bigl#1}{\bigr}{\,#3}{\bigr]}}
\newcommand{\SETDESCRBRACKETS}[3][|]%
     {\twincommandJN{\left[}{#2}{\left#1}{\right}{\,#3}{\right]}}

\newcommand{\olnot}[1]{\overline{#1}}

\newcommand{\synteq}{\doteq}

\newcommand{\nvar}{n}

\newcommand{\nclause}{m}

\newcommand{\clwidth}{k}

\newcommand{\randkcnfnclwrepl}[3][\clwidth]%
        {\ensuremath{\mathcal{F}^{#2, #3}_{#1}}}
\newcommand{\randkcnfnclwreplstd}%
        {\randkcnfnclwrepl{\clwidth}{\nvar}{\nclause}}

\newcommand{\complclassformat}[1]%
        {\textrm{\upshape{\textsf{#1}}}\xspace}
\newcommand{\cocomplclass}[1]%
        {\textrm{\upshape{\textsf{co#1}}}\xspace}

\newcommand{\DTIMEadviceclass}[2]%
    {\ensuremath{\complclassformat{DTIME}\bigl(#1\bigr)/{#2}}}

\newcommand{\PCPalph}[5]%
    {\ensuremath{\complclassformat{PCP}_{{#1},{#2}}[{#3}, {#4}, {#5}]}}
\newcommand{\PCP}[4]%
    {\ensuremath{\complclassformat{PCP}_{{#1},{#2}}[{#3}, {#4}]}}

\newcommand{\eqperiod}{\enspace .}
\newcommand{\eqcomma}{\enspace ,}
\renewcommand{\eqperiod}{\, .}
\renewcommand{\eqcomma}{\, ,}

\ifthenelse{\boolean{detectedToC}}{}
  { 
    }
\ifthenelse{\boolean{detectedToC}}{}{

}

\ifthenelse{\boolean{detectedLIPIcs} \or \boolean{detectedIJCAI}}
{\renewcommand{\st}{\errmessage{Please do not use st}}}
{\ifthenelse{\isundefined{\st}}
  {\newcommand{\st}{such that\xspace}}
  {}}      

\ifthenelse{\boolean{detectedIEEE}}{}{}

\ifthenelse
{\isundefined{\refeq}}
{\newcommand{\refeq}[1]{\eqref{#1}}}
{\renewcommand{\refeq}[1]{\eqref{#1}}}

\newcommand{\derives}{\vdash}

\ifthenelse{\boolean{detectedACMconf}}
{}
{}

\ifthenelse{\boolean{detectedACMconf}}
{
 }
{
 }

\newcommand{\SETSOFVARSORLIT}[2]%
        {\mathit{#1}\left({#2}\right)}
\newcommand{\setsofvarsorlit}[2]%
        {\mathit{#1}({#2})}
\newcommand{\Setsofvarsorlit}[2]%
        {\mathit{#1}\bigl({#2}\bigr)}

\newcommand{\restrict}[2]{{{#1}\!\!\upharpoonright_{#2}}}

\newcommand{\derivabbrev}[2]{\bigl( #1 \vdash #2 \bigr)}
\newcommand{\derivabbrevsmall}[2]{( #1 \vdash #2 )}
\newcommand{\derivabbrevcompact}[2]{\bigl( #1 \vdash #2 \bigr)}

\newcommand{\refutabbrevsmall}[1]{\derivabbrevsmall{#1}{\!\bot}}
\newcommand{\refutabbrevcompact}[1]{\derivabbrevcompact{#1}{\!\bot}}

\newcommand{\genericrefsmall}[3]%
    {{\mathit{#1}}_{#2}\refutabbrevsmall{#3}}
\newcommand{\genericrefcompact}[3]%
    {{\mathit{#1}}_{#2}\refutabbrevcompact{#3}}
\newcommand{\genericderiv}[4]%
    {{\mathit{#1}}_{#2}\derivabbrev{#3}{#4}}
\newcommand{\genericderivsmall}[4]%
    {{\mathit{#1}}_{#2}\derivabbrevsmall{#3}{#4}}
\newcommand{\genericderivcompact}[4]%
    {{\mathit{#1}}_{#2}\derivabbrevcompact{#3}{#4}}
\newcommand{\generictaut}[3]%
    {{\mathit{#1}}_{#2}\derivabbrev{}{#3}}
\newcommand{\generictautcompact}[3]%
    {{\mathit{#1}}_{#2}\derivabbrevcompact{}{#3}}
\newcommand{\generictautsmall}[3]%
    {{\mathit{#1}}_{#2}\derivabbrevsmall{}{#3}}

\newcommand{\formulaformat}[1]{\mathit{#1}}

\newcommand{\extendedversion}[1]{\widetilde{#1}}

\newcommand{\epopnot}[1]%
    {\extendedversion{\formulaformat{POP}}_{#1}}

\newcommand{\elopnot}[1]%
    {\extendedversion{\formulaformat{LOP}}_{#1}}

\newcommand{\ephpnot}[2]%
    {\vphantom{\extendedversion{\formulaformat{PHP}}}
      {\smash{\extendedversion{\formulaformat{PHP}}}
        \vphantom{\formulaformat{PHP}}}^{#1}_{#2}}

\newcommand{\efphpnot}[2]%
    {\vphantom{\extendedversion{\formulaformat{FPHP}}}
      {\smash{\extendedversion{\formulaformat{FPHP}}}
        \vphantom{\formulaformat{FPHP}}}^{#1}_{#2}}

\newcommand{\ontophpnot}[2]%
    {\formulaformat{Onto}\text{-}\formulaformat{PHP}^{#1}_{#2}}
\newcommand{\ontofphpnot}[2]%
    {\formulaformat{Onto}\text{-}\formulaformat{FPHP}^{#1}_{#2}}

\newcommand{\graphontophpnot}[1][G]%
    {\text{$\formulaformat{Onto}$-$\formulaformat{PHP}$}({#1})}

\newcommand{\perfectmatchingnot}[1][G]%
    {\formulaformat{PM}({#1})}
\newcommand{\core}{\mathcal{C}}
\newcommand{\derived}{\mathcal{D}}

\newcommand{\obj}{f}

\newcommand{\objorig}{\obj}
\newcommand{\objmodified}{\obj^\prime}

\newcommand{\newsubscript}{\mathit{new}}

\newcommand{\zeroone}{\mbox{$0$--$1$}\xspace}

\newcommand{\solver}[1]{\textsc{#1}\xspace}

\newcommand{\pb}{\solver{PB16}}
\newcommand{\miplib}{\solver{MIPLIB}}

\newcommand{\papilo}{\solver{PaPILO}}
\newcommand{\veripb}{\solver{Veri\-PB}}

\newcommand{\myorcidlink}[1]{\,\href{https://orcid.org/#1}{\raisebox{-0.45ex}{\includegraphics[width=1.8ex]{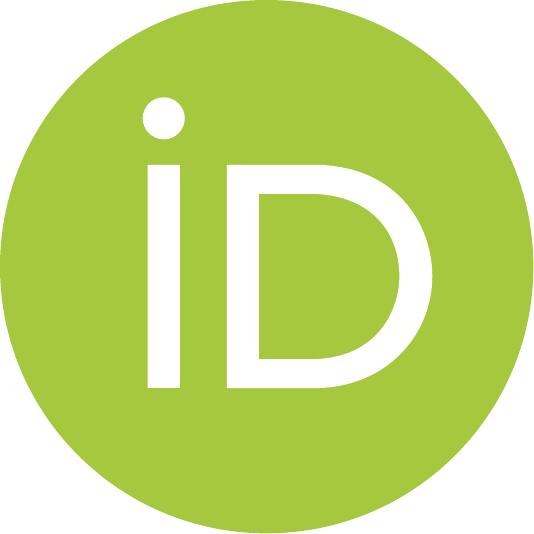}}}}

\definecolor{pbblue}{RGB}{68,119,170}
\definecolor{pbcyan}{RGB}{102,204,238}
\definecolor{pbgreen}{RGB}{34,136,51}
\definecolor{pbred}{RGB}{238,102,119}

\usepackage{todonotes,xspace}

\makeatletter
\def\orcidID#1{\href{http://orcid.org/#1}{\protect\raisebox{-1.25pt}{\protect\includegraphics{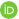}}}}
\makeatother

\begin{document}

\title{Certifying MIP-Based Presolve Reductions
  \\
  for
  \zeroone Integer Linear Programs
}
\author{
	Alexander Hoen\inst{1}\orcidID{0000-0003-1065-1651} \and
	Andy Oertel\inst{3,4}\orcidID{0000-0001-9783-6768} \and
	Ambros Gleixner\inst{1,2}\orcidID{0000-0003-0391-5903} \and 
	Jakob~Nordstr\"om\inst{3,4}\orcidID{0000-0002-2700-4285}}
\institute{%
  Zuse Institute Berlin, Takustr.~7, 14195 Berlin, Germany\\
  \email{hoen@zib.de}
  \and
  HTW Berlin, 10313 Berlin, Germany\\
  \email{gleixner@htw-berlin.de}
  \and 
  Lund University, Lund, Sweden\\
  \email{andy.oertel@cs.lth.se}\\
  \and
  University of Copenhagen, Copenhagen, Denmark \\
  \email{jn@di.ku.dk}  \\
}

\authorrunning{A. Hoen, A. Oertel, A. Gleixner, J. Nordstr\"om}

\maketitle

\begin{abstract}
It is
well known    
that
reformulating the original problem
can be
crucial for the performance of mixed-integer programming (MIP)
solvers.
To ensure correctness, all transformations must preserve the
feasibility status and optimal value of the problem, but there is
currently no established methodology to express and verify the
equivalence of two mixed-integer programs.
In this work, we take a first step in this direction by showing how
the correctness of MIP presolve reductions on \zeroone integer linear
programs can be certified by using (and suitably extending) the \veripb{}
tool for pseudo-Boolean proof logging.
Our experimental evaluation on both decision and optimization
instances demonstrates the computational viability of the approach and leads to
suggestions for
future revisions of the proof format
that will help to reduce the verbosity of the certificates and to
accelerate the certification and verification process further.

\end{abstract}
    

\section{Introduction}
\label{sec::introduction}


\emph{Boolean satisfiability solving (SAT)}
and
\emph{mixed-integer  programming (MIP)}
are two computational paradigms in which surprisingly mature and
powerful solvers have been developed over the last decades.
Today such solvers are routinely used to solve large-scale
problems in practice despite the fact that these problems are
\emph{NP}-hard.
Both SAT and MIP solvers typically start by trying to
simplify the input problem before feeding it to the main
solver algorithm, a process known as
\emph{presolving} in MIP
and
\emph{preprocessing} in SAT.
This can involve, e.g.,  fixing variables to values, strengthening
constraints, removing constraints, or adding new constraints to break
symmetries.
Such techniques are very important for SAT solver performance~%
\cite{BJK21PreprocessingPlusCrossref},
and for MIP solvers they often play a decisive role in whether a problem
instance can be solved or not, regardless of whether the solver uses
floating-point~\cite{presolving_achterberg}
or exact rational arithmetic~\cite{EiflerGleixner2022_Acomputational}.

The impressive performance gains for modern combinatorial solvers
come at the price of ever-increasing complexity, which makes these
tools very hard to debug. It is well documented that even
state-of-the-art solvers in many paradigms, not just SAT and MIP,
suffer from errors such as mistakenly claiming infeasibility or
optimality, or even returning ``solutions'' that are infeasible~%
\cite{%
  AGJMN18Metamorphic,%
  CKSW13Hybrid,%
  GSD19SolverCheck,%
  Klotz2014,%
  Steffy2011}.
During the last decade, the SAT community has dealt with this problem
in a remarkably successful way by requiring that solvers should use
\emph{proof logging}, i.e., produce machine-verifiable certificates
of correctness for their computations that can be verified by a
stand-alone proof checker. A number of proof formats have been
developed, such as
\solver{DRAT}~\cite{HHW13Trimming,HHW13Verifying,WHH14DRAT},
\solver{GRIT}~\cite{CMS17EfficientCertified},
and
\solver{LRAT}~\cite{CHHKS17EfficientCertified},
which are used to certify the whole solving process including
preprocessing.

Achieving something similar in a general MIP setting is much more challenging, amongst
others because of the presence of continuous and general integer variables, which may even
have unbounded domains.
For numerically exact MIP solvers~\cite{CookKochSteffyWolter2013,EiflerGleixner2022_Acomputational,EiflerGleixner2024_Safeandverified}
the proof format
\solver{VIPR}~\cite{CGS17Verifying}
has been introduced, but it currently only allows verification of feasibility-based reasoning, which must preserve all feasible solutions.
In particular, it does not support the verification
of dual presolving techniques that may exclude feasible solutions as long as one optimal solution remains. This means that while exact MIP solvers
could in principle generate a certificate for the main solving process,
such a certificate would only establish correctness under the
assumption that all the presolving steps were valid, as, e.g., in~\cite{EiflerGleixner2022_Acomputational}.
And, unfortunately, the proof logging techniques for SAT preprocessing
cannot be used to address this problem,
since they can only reason about
clausal constraints.

\subsection*{Our Contribution}

\noindent
In this work, we take a first step towards verification of the full
MIP solving process by demonstrating how pseudo-Boolean proof logging
with  \veripb can be used to produce certificates of correctness for a
wide range of MIP presolving techniques for \zeroone integer linear
programs (ILPs).
\veripb is quite a versatile tool in that it has previously been employed
for certification of, e.g.,  
advanced SAT solving techniques~\cite{BGMN23Dominance,GN21CertifyingParity},
SAT-based optimization (MaxSAT)~\cite{BBNOV23CertifiedCoreGuided,VWB22QMaxSATpb},
subgraph solving~%
\cite{GMMNPT20CertifyingSolvers,GMN20SubgraphIso},
and constraint programming~%
\cite{GMN22AuditableCP,MM23ProofLogging}.
However, to the best of our knowledge this is the first time the tool
has been used to prove the correctness of reformulations of
optimization problems, and this presents new challenges.
In particular, the proof system turns out not to be well suited for
problem reformulations with frequent changes to the objective
function, and therefore we introduce a new rule for objective function
updates.


Our computational experiments 
confirm that this approach to
certifying presolve reductions is computationally viable and
the overhead for certification aligns with what is known from the literature for certifying problem transformations in other contexts~\cite{GMNO22CertifiedCNFencodingPB}.
The analysis of the results reveals new insights into
performance bottlenecks, and these insights directly translate to possible revisions of the proof logging format that
would be valuable to address in order to decrease the size of the
generated proofs and speed up proof verification. 
%

We would like to note that, while our current methods are only applicable to \zeroone ILPs, this covers already a large and important class of MIPs.
In particular, there are applications where the exact and verified solution of \zeroone ILPs is highly relevant, see \cite{cip_phd,EiflerGleixnerPulaj2022,SahraouiBendottiAmbrosio2019} for some examples. 

The rest of this paper is organized as follows.
After presenting pseudo-Boolean proof logging and \veripb{} in
Sec.~\ref{sec::extension},
we demonstrate in
Sec.~\ref{sec::certifying}
how to produce \veripb certificates for
MIP presolving on \zeroone ILPs.
In
Sec.~\ref{sec::experiments}
we report results of an experimental evaluation,
and we conclude in
Sec.~\ref{sec::conclusion}
with a summary and discussion of future work.

\section{Pseudo-Boolean Proof Logging with \veripb}
\label{sec::extension}

We start by reviewing pseudo-Boolean reasoning
in Sec.~\ref{sec::pb-prelims}, and then explain
our extension to deal with objective function updates
in Sec.~\ref{subsubsec:obju}.

\subsection{Pseudo-Boolean Reasoning with  the Cutting Planes Method}
\label{sec::pb-prelims}

Our treatment of this material will by necessity be somewhat terse---we refer
the reader
to~\cite{BN21ProofCplxSATplusCrossref} for more information about the
cutting planes method and 
to~\cite{BGMN23Dominance,GMNO22CertifiedCNFencodingPB}
for detailed information about the  \veripb proof system and format.

We write~$x$ to denote a $\set{0,1}$-valued variable and
$\olnot{x}$ as a shorthand \mbox{for $1-x$},
and write~$\ell$ to denote such \emph{positive} and
\emph{negative literals}, respectively.
By a
\emph{pseudo-Boolean (PB) constraint}
we mean a \zeroone linear inequality
$\sum_j a_j \ell_j \geq b$,
where when convenient we can assume all literals~$\ell_j$ to refer to
distinct variables and all~$a_j$ and~$b$ to be 
non-negative
(so-called \emph{normalized form}).
A \emph{pseudo-Boolean formula} is 
just 
another name for a
\zeroone integer linear program.
For optimization problems we also have an
objective function $\obj = \sum_j c_j x_j$
that should be minimized (and $\obj$ can be negated to represent a maximization problem).

The foundation of \veripb is the
\emph{cutting planes} proof system~\cite{CCT87ComplexityCP}.
At the start of the proof, the set of 
\emph{core constraints}~$\core$
are initialized as the \zeroone linear inequalities in the problem
instance.
Any constraints derived as described below are placed in the
set of \emph{derived constraints}~$\derived$,
from where they can later be moved to~$\core$
(but not vice versa).
Loosely speaking, \veripb proofs maintain the invariant that
the optimal value of any solution to~$\core$ and to the original input
problem is the same.
New constraints can be derived from
$\core \cup  \derived$
by performing \emph{addition} of two constraints
or \emph{multiplication}
of a constraint by a positive integer,
and
\emph{literal axioms} $\ell \geq 0$
can be used at any time.
Additionally, for a constraint 
$\sum_{j} a_j \ell_j \geq b$
written in normalized form we can apply \emph{division} by a positive
integer~$d$ followed by rounding up to obtain
$\sum_{j} \ceiling{a_j / d} \ell_j \geq \ceiling{b / d}$,
and \emph{saturation} can be applied to yield
$\sum_j \min \set{a_j, b} \cdot \ell_j \geq b$.

For a PB constraint
$C \synteq \sum_{j} a_j \ell_j \geq b$
(where we use $\synteq$ to denote syntactic equality),
the negation of $C$
is
$
\neg C \synteq
\sum_{j} a_j \ell_j \leq b - 1$.
For a \emph{partial assignment}~$\rho$
mapping variables to~$\set{0,1}$,
we write
$\restrict{C}{\rho}$
for the \emph{restricted constraint} obtained by replacing variables
in~$C$ assigned by~$\rho$ by their values and  simplifying  the result.
We say that $C$
\emph{unit propagates $\ell$ under $\rho$} if $\restrict{C}{\rho}$
cannot be satisfied unless $\ell$ is assigned to~$1$.
If unit propagation on all constraints in
$\core \cup \derived \cup \set{\neg C}$
starting with the empty assignment
$\rho= \emptyset$,
and extending $\rho$ with new assignments as long as new literals propagate,
leads to contradiction in the form of a violated constraint,
then we say that
$C$ follows by \emph{reverse unit propagation (RUP)} from
$\core \cup \derived$.
Such
(efficiently verifiable)
RUP steps are allowed in \veripb proofs when it is convenient to
avoid writing out an explicit derivation of~$C$
from $\core \cup \derived$.
We will also write
$\restrict{C}{\omega}$
to denote the result of applying to~$C$ a
\emph{(partial) substitution}~$\omega$ which can remap variables to
other literals in addition to~$0$ and~$1$,
and we extend this notation to sets in the obvious way by taking unions.

In addition to the cutting planes rules, which can only derive
semantically implied constraints,  \veripb has a
\emph{redundance-based strengthening rule} that can derive a
non-implied constraint~$C$ as long as this does not change the feasibility
or optimal value of the problem.
Formally, $C$ can be derived from $\core \cup \derived$
using this rule by  exhibiting in the proof a
\emph{witness substitution}~$\omega$ together with subproofs
\begin{equation}
  \label{eq:redundance-rule}
  \core \cup \derived \cup \set{\neg C}
  \derives
  \restrict{(\core \cup \derived \cup \set{C})}{\omega}
  \cup
  \set{
    \obj \geq \restrict{\obj}{\omega}
  }
  \eqcomma
\end{equation}
of all constraints on the right-hand side from the premises on the
left-hand side using the derivation rules above.
Intuitively,
what \eqref{eq:redundance-rule} shows is
that if $\alpha$ is any assignment that
satisfies
$\core \cup \derived$ but violates~$C$, then
$\alpha \circ \omega$
satisfies
$\core \cup \derived \cup \set{C}$
and yields at least as good a value for the objective function~$\obj$.

During presolving, constraints in the input formula can be deleted or
replaced by other constraints, and the proof 
needs to
establish that such modifications are correct.
While deletions from the derived set~$\derived$ are always in order,
removing a constraint from the core set~$\core$ could potentially
introduce spurious solutions. Therefore, deleting a constraint~$C$
from~$\core$ can only be done by the 
\emph{checked deletion rule}, which requires 
to show that~$C$ could be rederived from $\core \setminus \set{C}$ by redundance-based strengthening
(see~\cite{BGMN23Dominance} for a more detailed explanation).

\subsection{A New Rule for Objective Function Updates}
\label{subsubsec:obju}

When variables are fixed or identified during the presolving process,
the objective function~$\objorig$
can be modified
to a
function~$\objmodified$%
.
This
modified objective~$\objmodified$ can then
be used in other presolver reasoning.
This scenario arises also in, e.g., MaxSAT solving, and can be dealt
with by deriving two PB constraints
$\objorig \geq \objmodified$
and
$\objmodified \geq \objorig$
in the proof%
, which encodes
that the old and new objective are 
equal~\cite{BBNOV23CertifiedCoreGuided}.
Whenever the solver argues in terms of~$\objmodified$, a
telescoping-sum argument with
$\objmodified = \objorig$ can be used to justify the same conclusion
in terms of the old objective.

However, if the presolver changes $\objorig$ to~$\objmodified$ and
then uses reasoning that needs to be certified by redundance-based strengthening,
then tricky problems can arise.
One of the required proof goals in~\eqref{eq:redundance-rule}
is that the witness~$\omega$ cannot worsen the objective.
If $\omega$ does not mention variables in~$\objmodified$, then this is
obvious to the presolver---$\omega$ has no effect on  the
objective---but if~$\omega$ assigns variables in the original
objective~$\objorig$, then one still needs to derive
$\objorig \geq \restrict{\objorig}{\omega}$
in the formal proof, which can be challenging.
While this can often be done by enlarging the witness~$\omega$ to
include earlier variable fixings and identifications, the extra
bookkeeping
required for this quickly becomes a major headache,
and results in the proof deviating further and further from
the actual presolver reasoning that the proof logging is meant to
certify.

For this reason, a better solution is to introduce a new 
\emph{objective function update rule}
that formally replaces~$\objorig$
by a new objective~$\objmodified$, so that all future reasoning about the
objective can focus on~$\objmodified$ and 
ignore~$\objorig$.
Such a rule needs to be designed with care, so that the optimal value
of the problem is preserved. 
Due to space constraints we cannot provide a formal proof here, but
recall that intuitively we maintain the invariant for the core
set~$\core$ that it has the same optimal value as the original
problem. In agreement with this, 
the formal requirement for
updating the objective from
$\objorig$ to~$\objmodified$
is to present in the proof log derivations of the two constraints
$\objorig \geq \objmodified$
and
$\objmodified \geq \objorig$
from the core set~$\core$ only.

\subsection{Examples for Proof Logging Syntax}

In order to make the concept of proof logging more concrete, we conclude this section by
providing, in Tab.~\ref{tab::veripb::explanation}, a few examples of how the derivation rules
explained above are encoded in \veripb{} syntax.
%
For space reasons, this list does not include examples of subproofs that may be necessary
for some derivations that cannot be proven automatically by \veripb.
Further details on practical aspects and implementation of pseudo-Boolean proof logging
can be found in the software repository of \veripb{}~\cite{veri_pb_dd7aa5a1}.

\begin{table}[t]
	\centering
	\scriptsize
	\caption{Examples of basic derivation rules in \veripb{} syntax. Here, (id) refers to the constraint ID assigned by \veripb.
	}
	\begin{tabular*}{\textwidth}{@{}l@{\;\;\extracolsep{\fill}}ll}
		\toprule
		\bf Rule & \bf Syntax & \bf Explanation\\[.5ex]
		\toprule
		\multirow{3}{*}{\parbox{2.5cm}{cutting planes in\\reverse Polish\\notation}} 
				&  \verb|pol x1 4 +| & add $x_1 \geq 0$ and (4)\\
				\cmidrule{2-3}
				&  \verb|pol 3 2 d|&  divides (3) by 2 \\
				\cmidrule{2-3}
				&  \verb|pol 1 2 * ~x1 +| &  multiplies (1) by 2 and adds $\olnot{x}_1 \geq 0$ \\
		\midrule
		\multirow{3}{*}{\parbox{2.5cm}{redundance-based\\strengthening\\}} 
				& \verb|red +1 x1 >= 1; x1 1|& verifies $x_1 \geq 1$ with $\omega = \{x_1\mapsto1\}$\\
				\cmidrule{2-3}
				& \verb|red +1 x1 +1 x2 >= 1; x1 x2 x2 x1|& verifies $x_1 + x_2\geq 1$ with\\ 
				&  & $\omega = \{x_1\mapsto x_2, x_2 \mapsto x_1\}$\\
		\midrule
		\multirow{1}{*}{\parbox{2.5cm}{RUP}} 
				& \verb|rup +1 x1 +1 x2 >= 1;|&verifies $x_1 + x_2\geq 1$ with RUP\\
		\midrule
		\multirow{1}{*}{\parbox{2.5cm}{move to core}} 
			& \verb|core id 3|&moves (3) to the core constraints \\
			\midrule
		\multirow{1}{*}{\parbox{2.5cm}{deletion from core}} 
			& \verb|delc 3|&deletes (3) from the core constraints \\
		\midrule
		\multirow{2}{*}{\parbox{2.5cm}{objective function\\update}} 
			& \verb|obju new +1 x1 +1 x2 1;|& defines $x_1 + x_2 +1$ as new objective\\
			\cmidrule{2-3}
			& \verb|obju diff +1 ~x1;|& adds $\olnot{x}_1$ to the objective\\
		\bottomrule
	\end{tabular*}
	\label{tab::veripb::explanation}
\end{table}

\newcommand{\ourparagraph}[1]{\medskip \noindent\textit{#1}.}

\section{Certifying Presolve Reductions}
\label{sec::certifying}

We now describe how feasibility- and optimality-based presolving
reductions can be certified by using \veripb proof logging enhanced
with the new objective function update rule described in Sec.~\ref{subsubsec:obju} above.
We distinguish
between
\emph{primal} and \emph{dual} reductions,
where
primal reductions strengthen the problem formulation by tightening the convex hull of the problem and preserve all feasible solutions,
and
dual reductions may additionally remove feasible solutions using optimality-based arguments.
More precisely, \emph{weak} dual reductions preserve all optimal solutions, but may remove suboptimal solutions.
\emph{Strong} dual reductions may remove also optimal solutions as long as at least one optimal solution is preserved in the reduced problem.
Our selection of methods is motivated by the recent
MIP solver
implementation described in~\cite{papilo}.
Before explaining the individual presolving techniques and their
certification, we introduce a few general techniques that are needed
for the certification of several presolving methods.

\subsection{General Techniques}
\label{subsec::general}


\ourparagraph{Substitution}
In order to reduce the number of variables, constraints, and non-zero coefficients in the constraints, many presolving techniques first try to identify an equality $E \synteq x_k = \sum_{j \not= k} \alpha_j x_j + \beta$ with $\alpha_j,\beta\in\Q$.
Subsequently, all occurrences of $x_k$ in the objective and constraints besides $E$ are substituted by the affine expression on the right-hand side and $x_k$ is removed from the problem.
The simplest case when $x_k$ is fixed to zero or one, i.e., when $\beta\in\{0,1\}$ and all $\alpha_j=0$, is straightforward to handle by deriving a new lower or upper bound on $x_k$.
During presolving, every fixed variable is removed from the model.
In the cases where some $\alpha_j\not= 0$, first the equation is expressed as a pair of constraints $E_{\geq} \land E_{\leq}$
%
and then the variable is removed by aggregation as follows.

\ourparagraph{Aggregation}
In order to substitute variables or reduce the number of non-zero coefficients, certain presolving techniques add a scaled equality $s \cdot E \synteq s\cdot E_{\geq}  \land s \cdot E_{\leq} $, $s \in \mathbb{Q}$, to a given constraint $D$. 
We call this an \emph{aggregation}.
Since \veripb{} certificates expect inequalities with integer coefficients, $s$ is split into two integer scaling factors $s_E, s_D\in \mathbb{Z}$ with $s = s_E/s_D$.
In the certificate, the aggregation is expressed as a newly derived constraint
	\begin{align*}
	D_{\newsubscript} \synteq
	\begin{cases}
		|s_E| \cdot E_{\geq} + |s_D| \cdot D & \;\text{if } \tfrac{s_D}{s_E} > 0\\
		|s_E| \cdot E_{\leq} + |s_D| \cdot D & \;\text{otherwise } \eqperiod
	\end{cases}
\end{align*}
Note that the presolving algorithm may decide to keep working with the constraint $(1/s_D) D_{\newsubscript}$ internally.
In this case, it must store the scaling factor $s_D$ in order to correctly translate between its own state and the state in the certificate; this happens in the implementation used in Sec.~\ref{sec::experiments}.


\ourparagraph{Checked Deletion}
The derivation of a new constraint $D_{\newsubscript}$ can render a previous constraint $D$ redundant.
A typical example is the case of substituting a variable above.
In a (pre)solver, the previous constraint is overwritten, and in order to keep the constraint database in the proof aligned with the solver, one may want to delete the previous constraint from the proof.
In order to check the deletion of $D$, a subproof is required that proves its redundancy.
In most cases, this subproof contains the ``inverted'' derivation of  $D_{\newsubscript}$.
As an example, consider an aggregation $D_{\newsubscript} \synteq D + E_{\leq}$ with an equality $E \synteq E_{\leq} \land E_{\geq}$.
In this case, the subproof for the checked deletion is $D_{\newsubscript} + E_\geq$.
Unless stated otherwise, the new constraints are moved to the core and redundant constraints are always removed by inverting the derivation of the constraint that replaces them.
%


\subsection{Primal Reductions}
\label{subsec::primal}

Primal reductions can be certified purely by implicational reasoning.

\ourparagraph{Bound Strengthening}
This preprocessor~\cite{FUGENSCHUH200569,Savelsbergh_Preprocessing_and_Probing_Techniques} tries to tighten the variable domains by iteratively applying well-known \emph{constraint propagation} to all variables in the linear constraints.
Each reduced variable domain is communicated to the affected constraints and may trigger further domain changes. 
This process is continued until no further domain reductions happen or the problem becomes infeasible due to empty domains.
Specifically, for an inequality constraint
\begin{align}
\label{eq:propagation_bound_orig}
\sum_{j \in N} a_{j} x_j \geq b
\end{align}
with $a_{k} \not= 0$, we first underestimate $a_{k}x_k$ via
%
$a_{k} x_k \geq b - \sum_{j \not= k} a_{j} x_j  \geq b - \sum_{j \not= k, a_{j}>0} a_{j}$. 
%
If $a_{k}>0$, this yields the lower bound
\newcommand{\MyCeiling}[1]{\Bigl \lceil #1 \Bigr \rceil}
\begin{align}
\label{eq:propagation_bound}
x_k \geq \MyCeiling{\big( b - \sum_{j\not= k, a_{j} > 0} a_{j} \big) / a_{k}} \eqcomma
\end{align}
and if $a_{k}<0$ we can obtain an analogous upper bound on $x_k$.

The bound change can be proven either by RUP, or more explicitly by stating the additions and division needed to form \eqref{eq:propagation_bound} from \eqref{eq:propagation_bound_orig} and the bound constraints.
We analyze the effect of both variants in Sec.~\ref{subsec::analysis-on-propagation}.

\ourparagraph{Parallel Rows}
Two constraints $C_j$ and $C_k$ are parallel if a scalar $\lambda\in \mathbb{R^+}$ exists with $\lambda(a_{j1},\ldots,a_{jn},b_j) = (a_{k1},\ldots,a_{kn},b_k)$.
Hence, one of these constraints is redundant and can be removed from the model~\cite{presolving_achterberg,GemanderChenWeningeretal}.
The subproof for deleting
the redundant rows must contain the remaining parallel row and $\lambda$ to prove the redundancy. 
For a fractional $\lambda$ the two constraint are scaled to ensure integer coefficients in the certificate.

\ourparagraph{Probing}
The general idea of \textit{probing}~\cite{cip_phd,Savelsbergh_Preprocessing_and_Probing_Techniques} is to tentatively fix a variable $x_j$ to 0 or 1 and then apply constraint propagation to the resulting model.
Suppose $x_k$ is an arbitrary variable with $k\not=i$, then we can learn fixings or implications in the following cases:
\begin{enumerate}
	\item If $x_j=0$ implies $x_k= 1$ and $x_j=1$ implies $x_k=0$ we can add the constraint $x_j = 1 - x_k $.
	Analogously, we can derive $x_k = x_j$ in the case that $x_j=0$ implies $x_k= 0$ and $x_j=1$ implies $x_k=1$.
	\item If $x_j =0$ propagates to infeasibility we can fix $x_j=1$.
	Analogously, if $x_j=1$ propagates to infeasibility we can fix $x_j=0$.
	\item If $x_j=0$ implies $x_k= 0$ and $x_j=1$ implies $x_k=0$ we can fix $x_k$ to 0.
	Analogously, $x_k$ can be fixed to 1 if $x_j=0$ implies $x_k= 1$ and $x_j=1$ implies $x_k=1$.
\end{enumerate} 
Cases 1 and 2 can be proven with RUP.
To prove correctness of fixing $x_k = 1$ in Case 3 we first derive two new constraints $x_k + x_j \geq 1$ and $ x_k - x_j \geq 0$ in the proof log by RUP.
Adding these two constraints leads to $x_k \geq 1$.
To prove $x_k=0$ we derive the constraints $x_k + x_j \leq 0$ and $ x_k - x_j \leq 0$ leading to $x_k = 0$.

\ourparagraph{Simple Probing}
On equalities with a special structure, a more simplified version of probing called \emph{simple~probing}~\cite[Sec.~3.6]{presolving_achterberg} can be applied.
Suppose 
the equation
\begin{align*}
  \sum_{j\in N} a_{j} x_j = b\,
  \text{ with }
  \sum_{j\in N} a_{j} = 2 \cdot b\,
  \text{ and }
	|a_{k}| = \sum_{j\in N, a_{j} >0} a_{j} - b
\end{align*}
holds for a variable $x_k$ with $a_{k} \not= 0$.  Let $\hat{N} = \set{p \in N \mid a_{p} \not= 0}$.
Under these conditions, $x_k = 1$ implies $x_p = 0$ and $x_k= 0$ implies $x_p = 1$ for all $p \in \hat{N}$ with $a_{p}>0$.
Further, $x_k = 1$ implies $x_p = 1$ and $x_k= 0$ implies $x_p = 0$ for all $p \in \hat{N}$ with $a_{p}<0$.
These implications can be expressed by the constraints
\begin{align}
	\label{eq::simple_probing}
	x_k = 1 - x_p &\text{ for all } p \in \hat{N} \textnormal{ with } a_{p} > 0 \eqcomma \\
		\label{eq::simple_probing2}
	x_k = x_p &\text{ for all } p \in \hat{N} \text{ with } a_{p} < 0 \eqperiod
\end{align}
The constraints \eqref{eq::simple_probing} and~\eqref{eq::simple_probing2} can be proven with RUP and used to substitute variables $x_p$ for all $p \in \hat{N}$ from the problem.

\ourparagraph{Sparsifying the Matrix}
The presolving technique sparsify~\cite{presolving_achterberg,Sparsify} tries to reduce the number of non-zero coefficients by adding (multiples of) equalities to other constraints using aggregations. This can be certified as described in Sec.~\ref{subsec::general}.

\ourparagraph{Coefficient Tightening}
The goal of this MIP presolving technique, which goes back to~\cite{Savelsbergh_Preprocessing_and_Probing_Techniques}, is to tighten the LP relaxation, i.e., the relaxation obtained when the integrality requirements are replaced by $x_j \in [0,1]$.
To this end, the coefficients of constraints are modified such that 
LP relaxation solutions are removed, but all integer feasible solutions are preserved.
%
Suppose we are given a constraint $\sum_{j \in N} a_{j} x_j \geq b$ with
$a_k \geq \varepsilon := a_k - b + \sum_{j \not= k, a_{j} < 0} a_{j} > 0$, 
then the constraint can be strengthened to
\begin{align*}
  (a_{k} - \varepsilon) x_k + \sum_{j \not= k} a_{j} x_j \geq b \eqperiod
\end{align*}
The case $a_k < 0$ is handled analogously.
This technique is also known as \emph{saturation} in the SAT community~\cite{CK05FastPseudoBoolean} and \veripb provides a dedicated saturation rule that can be used directly for proving the correctness of coefficient tightening.
The deletion of the original, weaker constraint can be proven automatically.

\ourparagraph{GCD-based Simplification}
This presolving technique from~\cite{Weninger2016} uses a divisibility argument to first eliminate variables from a constraint and then tighten its right-hand side.
Given $C \synteq \sum_{j \in N} a_{j} x_j \geq b$ with $|a_{1}|\geq\dots\geq |a_{n}|>0$.
We define the greatest common divisor $g_k = \gcd(a_1, \dots, a_k)$ as the largest value~$g$ such that $a_j/g \in \mathbb{Z}$ for all $j\in \{1,\dots,k\}$.
If for an index $k$ it holds that
\begin{align*}
b - g_k \cdot \CEILING{\frac{b}{g_k}} \geq \sum_{k< j\leq n, a_{j}>0} a_{j}
\quad\text{ and }\quad
b - g_k \cdot \CEILING{\frac{b}{g_k}} - g_{k} \leq \sum_{k< j\leq n, a_{j}<0} a_{j} \eqcomma
\end{align*}
then all $a_{k+1}, \dots, a_{n}$ can be set to 0. 
This first step can be certified as \textit{weakening}~\cite{weakening} and \veripb{} provides an out-of-the-box verification function for it.
Finally, $b$ can be rounded to $g_{k} \cdot \CEILING{b/g_{k}}$.
This rounding step can be certified by dividing $C$ with $g_k$ and then multiply it again with $g_{k}$.

\ourparagraph{Substituting Implied Free Variables}
A variable $x_j$ is called \emph{implied free} if its lower bound and its upper bound can be derived from the constraints.
For example, the constraints $x_1 - x_2 \geq 0$ and $x_2 \geq 0$ imply the lower bound $x_1\geq 0$.
%
If we have an implied free variable $x_j$ in an equality
$E \synteq a_{j} x_j + \sum_{k\not=j} a_{k} x_k = b$
with $a_{j} > 0$,
then we can remove $x_j$ from the problem by substituting it with 
$x_j = \big( b - \sum_{k\not=j} a_{k} x_k  \big) / a_{j}$, see~\cite{presolving_achterberg} for details.

To apply the substitution in the certificate we use aggregations to remove $x_j$ from all constraints and the objective function update to remove $x_j$ from the objective.
If coefficients $c_j/a_{j}$ or $a_{k}/a_{j}$ are non-integer,
then the resulting constraints are scaled as described in~Sec.~\ref{subsec::general}.
To prove the deletion of $E$, we 
derive two 
constraints by adding $x_j\geq 0$ and $1\geq x_j$ to $E$ each, which results in
\begin{align}
\label{eq::substitution_auxilary}
b \geq \sum_{k \not= j} a_{k} x_k
\;\land\;
\sum_{k \not= j} a_{k} x_k \geq b - a_{j} \eqperiod
\end{align}
Then the deletion of $E_{\geq}$ can be certified by a witness $\omega = \set{x_j \mapsto 1}$.
The constraint simplifies to~\eqref{eq::substitution_auxilary} and is therefore fulfilled.
Analogously, we use the witness $\omega = \set{x_j \mapsto 0}$ to certify the deletion of $E_{\leq}$.
Finally, to delete the constraints in~\eqref{eq::substitution_auxilary} we generate a subproof that
shows that negation of the auxiliary constraints in~\eqref{eq::substitution_auxilary} leads to $x_j\not\in \{0,1\}$.
This is a contradiction to the implied variable bounds $0\leq x_j \leq 1$.
Since these bounds are still present through the implying constraints, we can add these implying constraints to \eqref{eq::substitution_auxilary} in the subproof to arrive at a contradiction.

\ourparagraph{Singleton Variables}
It is well-known that variables that appear only in one inequality constraint or equality can be removed from the problem~\cite[Sec.~5.2]{presolving_achterberg}.
This can be certified by applying one of the following primal or dual strategies in this order:
First, try to apply duality-based fixing, see Sec.~\ref{subsec::dual};
second, an implied free singleton variable can be substituted as explained above;
otherwise, the singleton variable can be treated as a \emph{slack variable}:
substitute the variable in the objective, then relax the equality as in \eqref{eq::substitution_auxilary}, and delete the original constraint.

\subsection{Dual Reductions}
\label{subsec::dual}
Dual reductions remove solutions while preserving at least one optimal solution.
Hence, to prove the correctness of dual reductions we need to involve the redundance-based strengthening rule of \veripb{}.
For each derived constraint $C$ we only explain how to prove $\obj \geq \restrict{\obj}{\omega}$ (subject to the negation $\neg C$); the proof goals for $\restrict{C}{\omega}$ can be derived in a very similar fashion.

\ourparagraph{Duality-based Fixing}
This presolving step described in~\cite[Sec.~4.2]{presolving_achterberg} counts the \emph{down-} and \emph{up-lock} of a variable.
A down-lock on variable $x_j$ is a negative coefficient, an up-lock on variable $x_j$ is a positive coefficient (for $\geq$ constraints).
If $x_j$ has no down-locks and $c_j\leq 0$, it can be fixed to zero;
if $x_j$ has no up-locks and $c_j\geq 0$, it can be fixed to one.
%
These reductions can be certified with redundance-based strengthening using the witness $\omega = \set{x_j \mapsto v}$, where $v$ is the fixing value.
The proof goal for $\obj \geq \restrict{\obj}{\omega}$ is equivalent to
	$c_j x_j \geq c_j v$,
which is fulfilled by the conditions of duality-based fixing.

\ourparagraph{Dominated Variables}
A variable $x_j$ is said to \emph{dominate} another variable $x_k$~\cite{presolving_achterberg,Gamrath_ProgressInPresolving}, in notation $ x_j \succ x_k$, if
\begin{align}
c_j \leq c_k
\;\land\;
a_{ij} \geq a_{ik} \text{ for all } i \in \set{1,\dots,m}\eqcomma \label{eq::dom_constraints}
\end{align}
where $a_{ij}$ and $a_{ik}$ are the coefficients of variable $x_j$ and $x_k$, respectively, in the $i$-th constraint.
Variable $x_j$ is then favored over $x_k$ since $x_j$ contributes less to the objective function, but more to the feasibility of the constraints.
For every domination $x_j \succ x_k$, a constraint $C \synteq x_j \geq x_k$ can be introduced.
This constraint can be certified by redundance-based strengthening with the witness $\omega= \set{x_k \mapsto x_j, x_j \mapsto x_k}$.
The proof goal for $f \geq \restrict{f}{\omega}$ is equivalent to
\begin{align}
\label{eq::domcol_proofgoal}
c_j x_j + c_k x_k \geq c_j x_k + c_k x_j \eqperiod
\end{align}
The negated constraint $\neg C \synteq x_j < x_k$ leads to $x_k=1$ and $x_j=0$.
Substituting these values in~\eqref{eq::domcol_proofgoal} leads to $c_k \geq c_j$, which follows directly from Condition~\eqref{eq::dom_constraints}.

\ourparagraph{Dominated Variables Advanced}
For an implied free variable we can drop the variable bounds and pretend the variable is unbounded.
This allows for additional fixings in the following cases of dominated variables:
\begin{enumerate}[(a)]
	\item If the upper bound of $x_j$ is implied and $x_j \succ x_k$, then $x_k = 0$. \label{item:dom_advanced_1}
	\item If the lower bound of $x_k$ is implied and $x_j \succ x_k$, then $x_j = 1$. \label{item:dom_advanced_2}
	\item If the upper bound of $x_j$ is implied and $x_j \succ -x_k$, then $x_k = 1$. \label{item:dom_advanced_3}
	\item If the lower bound of $x_j$ is implied and $-x_j \succ x_k$, then $x_j = 0$. \label{item:dom_advanced_4}
\end{enumerate}
We use redundance-based strengthening with witness $\omega = \set{x_k \mapsto 0}$ to prove the correctness of~\ref{item:dom_advanced_1} as follows.
If the upper bound of $x_j$ is implied, this means there exists a constraint with  $a_{ij} < 0$ such that 
\newcommand{\MyFloor}[1]{\Bigl \lfloor #1 \Bigr \rfloor}
\begin{align}
x_j \leq \MyFloor{\frac{b_\ell - \sum_{\ell \not= j, a_{i\ell}>0} a_{i\ell}}{a_{ij}}} = 1 \eqperiod
\end{align}
Due to Condition~\eqref{eq::dom_constraints}, it must hold that $0 > a_{ij}\geq a_{ik}$, and the constraint 
$x_j + x_k \leq 1$
can be derived.
Hence, negating and propagating $C \synteq x_k = 0$ with RUP leads to contradiction, which proves the validity of $C$.
Case~\ref{item:dom_advanced_2} can be handled analogously using the witness $\omega = \set{x_k \mapsto 1}$.
To derive $C \synteq x_k=1$ in~\ref{item:dom_advanced_3} we use redundance-based strengthening with witness  $\omega = \set{ x_k \mapsto 1, x_j \mapsto 1}$.
Then, the proof goal for $f \geq \restrict{f}{\omega}$ is $c_j \cdot x_j + c_k \cdot x_k \geq c_j + c_k$.
After propagating $\neg C$, this becomes equivalent to $c_j \leq -c_k$, which is true by Condition~\eqref{eq::dom_constraints}. 
Case~\ref{item:dom_advanced_4} can be handled analogously using the witness $\omega = \set{x_k \mapsto 0,x_j \mapsto 0}$.

\subsection{Example}
\label{sec::example}

We conclude this section with an example
of a small certificate for the substitution of an implied free variable in Fig.~\ref{fig::example}, 
also available with a more detailed description at the software repository of \papilo~\cite{papilo_example}.
Consider the \zeroone ILP
\begin{align}
  \min\; x_1 + x_2\;
  \text{ s.t. }\;
	x_1 + x_2 - x_3 - x_4 &= 1
	\eqcomma
        \label{eq::example_equality}\\
	- x_1 + x_5 &\geq 0  
	\eqcomma
	\label{eq::example_ineq}
\end{align}
in which the lower bound of $x_1$ is implied by~\eqref{eq::example_equality} and the upper bound of $x_1$ is implied by~\eqref{eq::example_ineq}.
Hence, $x_1$ is implied free and we can use~\eqref{eq::example_equality} to substitute it.

In the left section of Fig.~\ref{fig::example} we first derive the two auxiliary constraints
\begin{align}
\label{eq::example_aux}
0 \leq x_2 - x_3 - x_4 \leq 1 \eqcomma
\end{align}
which receives the constraint IDs 4 and 5 and are moved to the core.
Note that the equality in \eqref{eq::example_equality} is split into two inequalities with IDs 1 and 2.
In the middle section, we first remove $x_1$ from \eqref{eq::example_ineq} by aggregation with~\eqref{eq::example_equality}, perform checked deletion,
then remove $x_1$ from the objective (automatically proven by \veripb).
Last, in the right section, we delete the equality in \eqref{eq::example_equality} used for the substitution and the auxiliary constraints in~\eqref{eq::example_aux} and arrive at the reformulated problem
\begin{align}
\min\; x_3 + x_4 +1 
\;\text{ s.t. }\; x_2 - x_3 - x_4 +x_5 \geq 1. 
\end{align}
From here, we could continue to derive $x_2=1$ by duality-based fixing, since $x_2$ has zero up-locks and objective coefficient zero.
This displays the importance of the objective update,
as without it
$x_2$ would still 
contribute to the objective with a positive coefficient, and this would prohibit duality-based fixing to~1.%

\begin{figure}[t]
  \centering
	\begin{boxedminipage}[t]{.32\linewidth}
		\begin{scriptsize}
\begin{verbatim}
* generates ID 4:
pol 1 ~x1 + ; 
core id 4
* generates ID 5:
pol 2 x1 + ;  
core id 5

\end{verbatim}\vspace*{-2.8ex}
		\end{scriptsize}
	\end{boxedminipage}
	\begin{boxedminipage}[t]{.32\linewidth}
		\begin{scriptsize}
\begin{verbatim}
* generates ID 6:
pol 3 1 + ;
core id 6       
delc 3 ;  ; begin
   pol 6 2 +
end
obju new +1 x3 +1 x4 1 ;
\end{verbatim}
		\end{scriptsize}
	\end{boxedminipage}
	\begin{boxedminipage}[t]{.32\linewidth}
		\begin{scriptsize}
\begin{verbatim}
delc 2 ; x1 -> 0
delc 1 ; x1 -> 1
delc 5
delc 4 ; ; begin
   pol 6 -1 +
end

\end{verbatim}\vspace*{-2.8ex}
		\end{scriptsize}
\end{boxedminipage}


	\caption{A \veripb{} certificate to substitute an implied free variable $x_1$.}
	\label{fig::example}
\end{figure}

\section{Computational Study}
\label{sec::experiments}

In this section we quantify the cost of \emph{certifying} presolve reductions in a state-of-the-art implementation for MIP-based presolve (Sec.~\ref{subsec:impact-of-proof-logging-on-papilo}) and the cost of \emph{verifying} the resulting certificates (Sec.~\ref{subsec::verifying-presolving}).
In Sec.~\ref{subsec::analysis-on-propagation} we analyze the impact of certifying constraint propagation by RUP or by an explicit cutting planes proof.

\subsection{Experimental Setup}

For generating the presolve certificates we use the solver-independent presolve library \papilo~\cite{papilo},
which provides a large set of MIP and LP techniques from the literature, described in Sec.~\ref{sec::certifying}.
Additionally, it accelerates the search for presolving reductions by parallelization,
encapsulating each reduction in a so-called transaction to avoid expensive synchronization~\cite{GGHpapilo}.
%
Logging the certificate, however, is performed sequentially while evaluating the transactions.

We base our experiments on models from the Pseudo-Boolean Competition 2016~\cite{PB16} including
1398~linear small integer decision and 532~linear small integer optimization instances of the competitions
PB10, PB11, PB12, PB15, and PB16 and 295~decision and 145~optimization instances from MIPLIB~2017~\cite{GleixnerHendelGamrathEtal2021} in the OPB translation~\cite{MIPLIB01}, excluding 10
large-scale instances\footnote{%
\solver{normalized-184},
\solver{normalized-pb-simp-nonunif},
\solver{a2864-99blp}, 
\solver{ivu06-big}, \solver{ivu59}, \solver{supportcase11}, \solver{a2864-99blp.0.s/u}, \solver{supportcase11.0.s/u}}
for which \papilo reaches the memory limit.
This yields a total of 671~optimization and 1681~decision instances.
We use \papilo~2.2.0~\cite{papilo_3b082d4}
running on 6 threads and \veripb~2.0~\cite{veri_pb_dd7aa5a1}.
The experiments are carried out on identical machines
with an 11th Gen Intel(R) Core(TM) i5-1145G7 @ 2.60 GHz CPU and 16 GB of memory and are assigned 14,000 MB of memory. 
The strict time limit for presolve plus certification and verification is three hours.
Times (reported in seconds) do not include the time for reading the instance file.
For all aggregations, we use the shifted geometric mean with a shift of 1~second.

\subsection{Overhead of Proof Logging}
\label{subsec:impact-of-proof-logging-on-papilo}

In the first experiment, we analyze the overhead of proof logging in \papilo.
The average results are summarized in Tab.~\ref{tab::papilo::overhead},
separately over decision (dec) and optimization (opt) instances for \pb and \miplib.  
Column ``relative'' indicates the average slow-down incurred by printing the certificate.

The relative overhead of proof logging is less than 6\% across all test sets.
\veripb supports two variants to change the objective function.
Either printing the entire 
objective (\verb|obju new|) or printing only the changes in the objective (\verb|obju diff|).
In our experiments, we only print the changes, since printing the entire objective for each change can lead to a large certificate and overhead, especially for instances with 
large and dense objective functions.
On the \pb instance \solver{normalized-datt256}, for example, \papilo finds $135\,206$ variable fixings.
Updating the entire objective function with $262\,144$ non-zeros for each of these variables leads to a huge certificate of about 138~GB and increases the time from 3.3~seconds (when printing only the changes) 
to 6625~seconds.\footnote{Certificate generated on Intel Xeon Gold 5122 @ 3.60GHz 96 GB with 50,000 MB of memory assigned.}

%
For 99\% of the instances, we can further observe that the \emph{overhead per applied reduction} is below $0.001\cdot 10^{-3}$~seconds over both test sets.
This means that the proof logging overhead is not only small on average, but 
also
small 
per applied reduction on the vast majority of instances.
These results show that the overhead scales well with the number of applied reductions and that proof logging remains viable even for instances with many transactions.
Here, under applied reductions we subsume all applied transactions and each variable fixing or row deletion in the first model clean-up phase.
During model clean-up, \papilo fixes variables and removes redundant constraints from the problem.
%
While \papilo technically does not count these reductions as full transactions found during the parallel presolve phase, their certification can incur the same overhead.


\begin{table}[t]
	\centering
	\small
	\caption{Runtime comparison of \papilo with and without proof logging.}
	\begin{tabular*}{\textwidth}{@{}l@{\;\;\extracolsep{\fill}}rrrr}
		\toprule
		test set & size &  default [s] & w/proof log [s] & relative \\
		\midrule
		\pb-dec & 1397 & 0.06 & 0.06 &  1.00\\[-0.5ex]
		\miplib-dec & 291 & 0.42 & 0.43 &  1.02\\[-0.5ex]
		\pb-opt & 531 & 0.65 & 0.66 &  1.02\\[-0.5ex]
		\miplib-opt & 142 & 0.33 & 0.35 &  1.06\\[-0.5ex]
		\bottomrule
	\end{tabular*}
	\label{tab::papilo::overhead}
\end{table}

\begin{table}[t]
	\centering
	\small
	\caption{
		Time to verify the certificates%
		.
		\veripb{} timeouts are treated with PAR2.
		}
	\begin{tabular*}{\textwidth}{@{}l@{\;\;\extracolsep{\fill}}rrrrrrr}
		\toprule
		&           &           & \multicolumn{2}{c}{\papilo time [s]} & \multicolumn{1}{c}{\veripb} & \multicolumn{2}{c}{relative time w.r.t.} \\
		\cmidrule{4-5}
		\cmidrule{7-8}
		test set & size &  verified  & default  & w/proof log& \multicolumn{1}{c}{time [s]} &  default &  w/proof log  \\
		\midrule
		\pb-dec & 1397 & 1397 & 0.06 & 0.06 & 0.88 & 14.67 & 14.67\\[-0.5ex]
		\miplib-dec & 291 & 267  & 0.42  & 0.43 & 9.64 & 22.85 & 22.42\\[-0.5ex]
		\pb-opt & 531 & 520 & 0.65 & 0.66 & 10.44 & 16.06 & 15.82 \\[-0.5ex]
		\miplib-opt & 142 & 139 & 0.33 & 0.35 & 5.25 & 15.91 & 15.00\\[-0.5ex]
		\bottomrule
	\end{tabular*}
	\label{tab::papilo_vs_veripb}
\end{table}

\subsection{Verification Performance on Presolve Certificates}
\label{subsec::verifying-presolving}

In this section, we analyze the time to verify the certificates generated by \papilo.
The results
are summarized in Tab.~\ref{tab::papilo_vs_veripb}.
The ``verified'' column lists the number of instances verified within 3~hours.
\veripb timeouts are counted as twice the time limit, i.e.,
PAR2 score.
Similar to Tab.~\ref{tab::papilo::overhead}, the ``relative'' columns report the relative overhead of \veripb{} runtime compared to \papilo.

First note that all certificates are
verified by \veripb{} (partially on the 38~instances where \veripb{} times out).
%
On average, it takes between $14.7$ and $22.4$ times as much time to verify the certificates than to produce them.
Nevertheless, some instances take a longer than average time to verify.
Over all test sets, 25\% of the instances have an overhead of at least a factor of $193$, see also Fig.~\ref{fig::experiment::PaPILOVeriPB}.

To put this result into context, note that presolving amounts more to a transformation than to a (partial) solution of the problem.
Each reduction has to be certified and verified while a purely solution-targeted algorithm may be able to skip certifying of a larger part of the findings that are not form a part of the final proof of optimality.
Hence, it makes sense to compare the performance of \veripb on presolve certificates to the overhead for, e.g., for verifying CNF translations~\cite{GMNO22CertifiedCNFencodingPB}.
For this study, a similar performance overhead is reported as in Fig.~\ref{fig::experiment::PaPILOVeriPB}.

\begin{figure}[t]
	\centering
	\begin{minipage}[t]{.40\linewidth}
		\begin{tikzpicture}[scale=0.6, transform shape]
		\begin{axis}[
		scatter/classes={blue={mark=x,draw=pbcyan},orange={mark=x,draw=gray},red={mark=x,draw=gray},green={mark=+,draw=pbgreen},purple={mark=x,draw=gray},olive={mark=x,draw=blue},peer={mark=x,draw=white}},
		ylabel=\veripb (time in seconds),
		xlabel=\papilo (time in seconds),
		xmode=log,
		ymode=log,
		grid=both,
		grid style={line width=.1pt, draw=gray!10},
		major grid style={line width=.2pt,draw=gray!50},
		minor tick num=5,
		xtick={-0.01,0.01,0.1,1,10,100,1000},
		ytick={-0.01,0.01,0.1,1,10,100,1000,10000},
		]
		\draw[domain=-100:100, smooth, variable=\x, black] plot ({\x}, {\x });
		\draw[domain=-100:100, smooth, variable=\x, black] plot ({\x}, {\x + 10});
		\addplot[scatter,only marks, scatter src=explicit symbolic]
		table[meta=label] {
			x y label
			1000 1000 peer
			-0.01 -0.01 peer
			3.798 96.03 green
			2.867 161.12 green
			2.414 146.46 green
			3.924 186.22 green
			1.709 188.82999999999998 green
			2.165 207.76999999999998 green
			5.941 133.74 green
			3.539 124.39 green
			2.576 119.23 green
			1.364 94.25 green
			2.092 119.78 green
			2.35 110.08 green
			0.825 158.07999999999998 green
			1.014 143.22 green
			2.098 124.86000000000001 green
			2.283 219.14999999999998 green
			2.192 136.81 green
			1.751 188.59 green
			1.671 106.8 green
			4.193 166.22 green
			0.879 133.7 green
			4.05 179.98 green
			6.445 178.19 green
			7.607 157.57 green
			0.912 168.62 green
			1.982 126.78999999999999 green
			2.459 100.95 green
			2.501 112.58 green
			3.334 133.74 green
			0.943 109.07000000000001 green
			3.539 135.28 green
			5.216 153.17000000000002 green
			2.684 117.33999999999999 green
			3.619 157.24 green
			4.549 138.55 green
			2.917 165.89 green
			1.497 167.38 green
			5.187 169.34 green
			1.925 151.36 green
			2.091 104.17999999999999 green
			2.444 182.03 green
			1.782 167.82 green
			6.105 170.35999999999999 green
			3.262 178.58 green
			1.505 120.88 green
			1.594 123.57 green
			2.32 189.91 green
			1.008 151.23 green
			2.445 119.26 green
			0.994 188.18 green
			1.121 175.36 green
			4.42 185.89999999999998 green
			7.15 152.85 green
			2.556 204.51000000000002 green
			7.726 178.32000000000002 green
			3.653 37.93 green
			1.025 34.65 green
			0.942 37.559999999999995 green
			1.587 25.709999999999997 green
			1.214 41.82 green
			1.316 29.52 green
			1.571 25.94 green
			0.994 36.85 green
			1.502 23.939999999999998 green
			0.943 37.41 green
			1.76 41.18 green
			1.453 28.26 green
			2.888 34.12 green
			2.447 37.0 green
			1.152 27.61 green
			0.801 30.62 green
			1.442 42.98 green
			2.398 22.19 green
			1.011 33.24 green
			1.345 57.37 green
			1.577 25.44 green
			1.725 37.48 green
			1.435 29.9 green
			2.843 78.36999999999999 green
			2.343 101.25 green
			3.279 52.51 green
			0.682 68.25 green
			2.577 89.07 green
			2.48 76.21000000000001 green
			2.047 82.0 green
			5.863 96.66 green
			2.763 96.31 green
			1.138 76.97999999999999 green
			2.553 98.92999999999999 green
			3.998 77.61 green
			1.024 91.56 green
			1.936 78.01 green
			3.86 81.03999999999999 green
			2.835 97.23 green
			0.665 67.56 green
			0.824 86.77 green
			2.013 121.46000000000001 green
			1.353 107.64 green
			1.938 55.279999999999994 green
			3.638 78.92999999999999 green
			1.499 73.98 green
			3.918 76.21000000000001 green
			3.773 83.96 green
			1.616 81.33 green
			1.237 74.85 green
			3.034 77.41 green
			1.99 59.589999999999996 green
			1.251 78.85 green
			1.516 80.19 green
			3.571 71.0 green
			0.638 9.03 green
			0.966 9.709999999999999 green
			0.559 10.42 green
			0.287 8.76 green
			0.471 10.809999999999999 green
			0.566 6.0 green
			0.836 6.630000000000001 green
			0.673 8.959999999999999 green
			1.109 9.6 green
			1.114 10.2 green
			0.693 8.479999999999999 green
			0.458 8.17 green
			0.422 6.640000000000001 green
			0.227 9.219999999999999 green
			2.011 11.82 green
			1.102 6.21 green
			0.555 8.19 green
			0.649 12.53 green
			0.599 10.629999999999999 green
			1.868 11.76 green
			0.837 8.719999999999999 green
			0.724 6.7 green
			1.012 10.25 green
			1.098 11.66 green
			1.028 6.94 green
			0.754 5.0 green
			0.735 6.46 green
			0.771 8.52 green
			0.744 8.24 green
			0.007 1.24 green
			0.007 1.26 green
			0.005 1.0 green
			0.067 6.92 green
			0.058 6.6499999999999995 green
			0.055 7.23 green
			0.064 7.819999999999999 green
			0.056 6.64 green
			0.058 7.7299999999999995 green
			0.064 5.16 green
			0.067 7.569999999999999 green
			0.054 7.66 green
			0.058 5.930000000000001 green
			0.07 6.79 green
			0.081 7.06 green
			0.063 6.75 green
			0.061 6.46 green
			0.013 1.42 green
			0.063 7.2 green
			0.063 6.9399999999999995 green
			0.051 6.1 green
			0.059 7.3 green
			0.055 6.21 green
			0.056 7.09 green
			0.052 7.12 green
			0.062 7.3999999999999995 green
			0.055 6.279999999999999 green
			0.06 7.31 green
			0.067 8.940000000000001 green
			0.067 9.870000000000001 green
			0.07 8.530000000000001 green
			0.053 4.510000000000001 green
			0.018 1.82 green
			0.057 8.100000000000001 green
			0.058 8.33 green
			0.017 1.95 green
			0.024 2.48 green
			0.07 8.42 green
			0.038 10.940000000000001 green
			0.06 8.99 green
			0.014 1.4 green
			0.062 8.350000000000001 green
			0.061 8.47 green
			0.017 2.54 green
			0.211 61.25 green
			0.832 638.84 green
			0.033 4.04 green
			49.232 421.58 green
			0.749 26.01 green
			10.665 13.219999999999999 green
			0.068 102.11 green
			1.058 1.7200000000000002 green
			0.01 1.83 green
			0.064 311.11 green
			0.008 2.9000000000000004 green
			0.018 1.55 green
			0.178 2519.28 green
			0.013 5.08 green
			0.009 1.54 green
			3.355 105.99000000000001 green
			0.254 21.23 green
			2.355 1.3800000000000008 green
			32.01 4605.26 green
			30.612 4658.59 green
			32.981 543.4100000000001 green
			25.358 9715.19 green
			30.69 4581.91 green
			26.286 9755.460000000001 green
			29.759 4663.110000000001 green
			24.339 9726.05 green
			26.496 9668.22 green
			29.875 4618.59 green
			30.033 4597.33 green
			1.756 9667.83 green
			25.997 9705.29 green
			29.926 4653.5199999999995 green
			1.936 4710.24 green
			25.187 9701.35 green
			31.403 4615.2 green
			25.185 9715.85 green
			30.901 4610.400000000001 green
			0.243 928.01 green
			0.252 448.49 green
			0.151 316.92 green
			0.151 314.71 green
			0.173 131.7 green
			0.175 124.66 green
			0.181 27.380000000000003 green
			0.255 450.98 green
			0.137 129.9 green
			0.205 1445.27 green
			0.223 929.81 green
			0.242 920.89 green
			0.151 316.33 green
			0.141 317.77 green
			0.229 922.41 green
			0.171 128.5 green
			0.171 125.33 green
			0.257 446.7 green
			0.204 1441.7 green
			0.18 74.34 green
			0.203 935.9699999999999 green
			0.163 2902.8799999999997 green
			0.182 75.11 green
			0.179 27.37 green
			0.139 322.59 green
			0.624 2887.3399999999997 green
			0.21 1443.99 green
			0.18 27.39 green
			0.174 131.94 green
			0.244 925.0999999999999 green
			0.701 248.27 green
			0.416 222.92999999999998 green
			0.691 254.24 green
			0.412 214.42 green
			1.28 276.47 green
			1.833 263.34999999999997 green
			0.493 108.72999999999999 green
			1.382 219.58 green
			0.468 243.64 green
			0.404 191.67000000000002 green
			1.302 255.67 green
			0.392 96.89 green
			0.42 224.53 green
			0.469 106.39 green
			0.696 257.89000000000004 green
			0.409 218.79 green
			0.805 240.7 green
			1.036 256.21999999999997 green
			0.369 99.41 green
			0.482 101.38 green
			0.091 2.87 green
			0.136 5.949999999999999 green
			0.043 1.05 green
			0.147 7.37 green
			0.055 16.84 blue
			0.034 20.64 blue
			0.055 15.21 blue
			0.019 9.46 blue
			0.022 11.24 blue
			0.034 21.169999999999998 blue
			0.036 17.68 blue
			0.036 20.27 blue
			0.014 5.5200000000000005 blue
			0.022 10.950000000000001 blue
			0.016 7.1899999999999995 blue
			0.042 26.9 blue
			0.015 5.930000000000001 blue
			0.088 25.87 blue
			0.034 19.93 blue
			0.016 7.5200000000000005 blue
			0.424 37.23 blue
			0.015 6.24 blue
			0.038 21.799999999999997 blue
			0.035 22.07 blue
			0.04 21.869999999999997 blue
			0.036 21.61 blue
			0.037 24.47 blue
			0.085 13.72 blue
			0.037 19.05 blue
			0.043 30.209999999999997 blue
			0.029 17.38 blue
			0.011 3.86 blue
			0.087 19.91 blue
			0.031 17.43 blue
			0.032 19.12 blue
			0.015 6.59 blue
			0.038 20.619999999999997 blue
			0.022 11.360000000000001 blue
			0.05 36.59 blue
			0.038 24.18 blue
			0.088 70.38 blue
			0.068 37.190000000000005 blue
			0.074 61.349999999999994 blue
			0.077 68.98 blue
			0.029 16.12 blue
			0.036 22.45 blue
			0.036 23.02 blue
			0.075 64.7 blue
			0.032 19.17 blue
			0.023 10.930000000000001 blue
			0.035 20.32 blue
			0.101 70.46000000000001 blue
			0.075 59.169999999999995 blue
			0.082 58.72 blue
			0.032 18.35 blue
			0.066 41.99 blue
			0.076 31.71 blue
			0.025 13.32 blue
			0.045 29.009999999999998 blue
			0.049 33.169999999999995 blue
			0.073 61.599999999999994 blue
			1.65 82.35000000000001 blue
			0.054 32.79 blue
			0.066 37.45 blue
			0.084 53.589999999999996 blue
			0.057 50.15 blue
			0.073 40.540000000000006 blue
			0.043 27.61 blue
			2.094 105.0 blue
			0.836 127.13 blue
			0.09 69.25 blue
			0.041 27.06 blue
			0.246 113.14 blue
			4.235 105.53 blue
			0.248 124.32000000000001 blue
			0.058 46.13 blue
			0.437 121.14999999999999 blue
			0.102 100.38 blue
			0.151 110.06 blue
			0.061 48.43 blue
			0.121 109.47999999999999 blue
			1.292 167.09 blue
			0.085 65.25 blue
			0.348 129.56 blue
			0.082 69.83000000000001 blue
			0.132 89.77 blue
			0.257 106.18 blue
			0.138 99.08 blue
			0.056 42.33 blue
			0.466 114.82000000000001 blue
			0.74 87.32000000000001 blue
			0.486 101.98 blue
			0.077 64.49000000000001 blue
			1.103 134.4 blue
			0.092 88.26 blue
			0.237 105.91 blue
			0.691 146.03 blue
			0.094 74.98 blue
			0.264 87.21 blue
			0.111 66.55000000000001 blue
			0.75 72.98 blue
			0.237 82.71 blue
			1.669 96.58 blue
			0.806 93.65 blue
			0.124 99.25 blue
			6.033 175.36 blue
			0.081 66.77000000000001 blue
			0.418 139.12 blue
			1.596 160.69000000000003 blue
			0.065 47.470000000000006 blue
			0.406 119.03 blue
			0.016 5.63 blue
			0.03 5.01 blue
			0.013 4.79 blue
			0.004 1.11 blue
			0.056 7.680000000000001 blue
			0.013 4.76 blue
			0.016 4.44 blue
			0.017 6.46 blue
			0.018 7.15 blue
			0.092 8.93 blue
			0.067 8.280000000000001 blue
			0.021 6.140000000000001 blue
			0.018 7.38 blue
			0.019 7.44 blue
			0.005 1.46 blue
			0.048 6.69 blue
			0.039 6.41 blue
			0.032 4.17 blue
			0.019 8.44 blue
			0.012 4.07 blue
			0.026 4.4 blue
			0.017 6.86 blue
			0.036 6.57 blue
			0.014 4.78 blue
			0.035 9.66 blue
			0.014 4.78 blue
			0.007 2.35 blue
			0.008 2.8400000000000003 blue
			0.012 5.180000000000001 blue
			0.014 4.8500000000000005 blue
			1.841 292.87 blue
			1.651 2.79 blue
			2.039 6.130000000000001 blue
			2.178 3.67 blue
			0.309 2.83 blue
			0.226 2.37 blue
			11.604 24.059999999999995 blue
			20.692 45.33 blue
			0.198 2.49 blue
			1.633 5.9399999999999995 blue
			1.869 3.2 blue
			2.17 3.58 blue
			1.867 3.19 blue
			547.066 827.75 blue
			1.916 3.1900000000000004 blue
			10.491 22.140000000000004 blue
			3.217 59.82 blue
			0.258 5.79 blue
			1.8 12.02 blue
			0.457 101.34 blue
			0.594 96.32000000000001 blue
			0.456 100.9 blue
			8.903 18.67 blue
			0.235 2.87 blue
			1.625 2.84 blue
			2.166 47.92 blue
			7.993 16.26 blue
			0.449 101.27000000000001 blue
			1.877 3.13 blue
			1.873 3.1399999999999997 blue
			1.574 7.300000000000001 blue
			0.194 2.0999999999999996 blue
			382.464 539.67 blue
			0.24 2.83 blue
			0.142 2.25 blue
			0.215 2.73 blue
			0.186 2.58 blue
			0.137 2.25 blue
			0.161 2.41 blue
			0.209 2.6 blue
			0.242 2.88 blue
			0.184 2.58 blue
			0.171 2.48 blue
			0.136 2.25 blue
			0.218 2.74 blue
			0.142 2.33 blue
			0.128 2.18 blue
			0.24 2.7399999999999998 blue
			0.154 2.41 blue
			0.191 2.5900000000000003 blue
			0.156 2.37 blue
			0.167 2.45 blue
			0.204 2.7 blue
			0.237 2.6900000000000004 blue
			0.139 2.23 blue
			0.151 2.39 blue
			0.212 2.72 blue
			0.132 2.28 blue
			0.196 2.5900000000000003 blue
			0.15 2.3 blue
			0.192 2.5100000000000002 blue
			0.206 2.7 blue
			0.19 2.4000000000000004 blue
			0.156 10.66 blue
			0.046 1.73 blue
			0.152 8.83 blue
			0.093 4.27 blue
			0.206 65.59 blue
			0.11 75.11 blue
			0.219 66.17999999999999 blue
			0.004 1.42 blue
			0.012 1.82 blue
			0.013 2.0500000000000003 blue
			0.012 1.92 blue
			0.012 2.0500000000000003 blue
			0.011 1.67 blue
			0.012 1.86 blue
			0.012 2.19 blue
			0.012 2.1300000000000003 blue
			0.011 1.84 blue
			0.011 1.58 blue
			0.012 2.22 blue
			0.012 1.65 blue
			0.011 2.18 blue
			0.012 2.1500000000000004 blue
			0.011 2.0500000000000003 blue
			0.012 1.66 blue
			0.011 1.79 blue
			0.012 1.86 blue
			0.011 1.8 blue
			0.011 1.81 blue
			0.011 2.14 blue
			0.011 1.9999999999999998 blue
			0.011 1.57 blue
			0.012 2.0 blue
			0.011 1.72 blue
			0.011 1.6199999999999999 blue
			0.012 1.8699999999999999 blue
			0.013 2.22 blue
			0.011 1.91 blue
			0.011 1.89 blue
			0.012 1.75 blue
			0.011 1.58 blue
			0.012 1.84 blue
			0.012 1.93 blue
			0.012 1.99 blue
			0.012 1.71 blue
			0.011 1.82 blue
			0.012 1.73 blue
			0.013 2.22 blue
			0.011 1.82 blue
			0.011 1.64 blue
			0.011 1.67 blue
			0.011 1.74 blue
			0.011 2.12 blue
			0.011 1.9 blue
			0.01 1.57 blue
			0.011 1.71 blue
			0.011 1.97 blue
			0.012 1.81 blue
			0.012 1.98 blue
			0.012 1.58 blue
			0.012 2.2300000000000004 blue
			0.012 2.0300000000000002 blue
			0.011 1.6199999999999999 blue
			0.012 1.84 blue
			0.011 1.59 blue
			0.011 1.64 blue
			0.012 1.76 blue
			0.012 2.2 blue
			0.012 2.22 blue
			0.011 1.9 blue
			0.011 1.9999999999999998 blue
			0.011 1.72 blue
			0.011 1.8699999999999999 blue
			0.011 1.98 blue
			0.01 1.71 blue
			0.012 1.92 blue
			0.011 1.91 blue
			0.012 1.84 blue
			0.013 1.71 blue
			0.012 2.2 blue
			0.011 1.81 blue
			0.012 2.0300000000000002 blue
			0.011 1.85 blue
			0.012 1.73 blue
			0.011 1.77 blue
			0.013 2.18 blue
			0.012 1.9999999999999998 blue
			0.011 2.06 blue
			0.01 1.67 blue
			0.011 1.6199999999999999 blue
			0.011 2.1900000000000004 blue
			0.011 1.9 blue
			0.012 2.04 blue
			0.011 1.88 blue
			0.011 1.89 blue
			0.012 2.15 blue
			0.011 1.73 blue
			0.011 1.93 blue
			0.012 1.65 blue
			0.011 1.94 blue
			0.012 2.16 blue
			0.012 1.93 blue
			0.011 1.61 blue
			0.011 1.6199999999999999 blue
			0.01 1.72 blue
			0.012 1.95 blue
			0.011 1.61 blue
			0.011 1.73 blue
			0.013 1.9899999999999998 blue
			0.012 1.98 blue
			0.011 1.63 blue
			0.012 1.67 blue
			0.01 1.69 blue
			0.012 2.03 blue
			0.012 1.6199999999999999 blue
			0.012 2.22 blue
			0.012 1.6 blue
			0.012 1.74 blue
			0.011 2.0100000000000002 blue
			0.011 1.9999999999999998 blue
			0.012 2.0300000000000002 blue
			0.012 1.71 blue
			0.011 1.81 blue
			0.011 1.8 blue
			0.011 1.88 blue
			0.011 2.0300000000000002 blue
			0.011 1.6 blue
			0.011 2.04 blue
			0.013 1.8599999999999999 blue
			0.012 1.9 blue
			0.012 1.9 blue
			0.011 1.72 blue
			0.013 1.57 blue
			0.012 1.85 blue
			0.011 1.6 blue
			0.012 2.1 blue
			0.012 1.85 blue
			0.012 2.06 blue
			0.011 1.88 blue
			0.012 2.02 blue
			0.012 2.24 blue
			0.011 1.6199999999999999 blue
			0.011 1.67 blue
			0.012 1.65 blue
			0.011 2.02 blue
			0.012 1.61 blue
			0.011 2.12 blue
			0.011 1.88 blue
			0.013 1.84 blue
			0.011 2.04 blue
			0.011 1.67 blue
			0.012 1.77 blue
			0.012 2.06 blue
			0.012 2.05 blue
			0.011 1.89 blue
			0.012 1.94 blue
			0.011 1.74 blue
			0.011 2.1700000000000004 blue
			0.013 2.2100000000000004 blue
			0.011 1.64 blue
			0.011 1.93 blue
			0.011 1.8 blue
			0.012 2.16 blue
			0.011 1.72 blue
			0.01 1.75 blue
			0.012 1.92 blue
			0.012 1.82 blue
			0.013 1.65 blue
			0.011 1.82 blue
			0.011 1.83 blue
			0.012 2.16 blue
			0.011 2.12 blue
			0.011 1.83 blue
			0.011 2.02 blue
			0.011 1.63 blue
			0.011 1.6199999999999999 blue
			0.011 1.6 blue
			0.012 2.17 blue
			0.012 1.76 blue
			0.011 1.9 blue
			0.011 1.57 blue
			0.012 2.22 blue
			0.012 2.16 blue
			0.012 1.92 blue
			0.011 2.0100000000000002 blue
			0.012 2.1900000000000004 blue
			0.011 1.83 blue
			0.013 1.75 blue
			0.011 1.68 blue
			0.044 12.760000000000002 blue
			0.017 5.37 blue
			0.009 1.3699999999999999 blue
			0.011 1.52 blue
			0.011 1.7 blue
			0.044 12.629999999999999 blue
			0.021 5.91 blue
			0.021 6.220000000000001 blue
			0.039 11.58 blue
		};
		\end{axis}
		\end{tikzpicture}
	\end{minipage}
	\begin{minipage}[t]{.05\linewidth}
		\ 
	\end{minipage}
	\begin{minipage}[t]{.40\linewidth}
		\begin{tikzpicture}[scale=0.6, transform shape]
		\begin{axis}[
		scatter/classes={blue={mark=x,draw=pbcyan},orange={mark=x,draw=gray},red={mark=x,draw=gray},green={mark=+,draw=pbgreen},purple={mark=x,draw=gray},olive={mark=x,draw=blue},peer={mark=x,draw=white}},
		ylabel=\veripb (time in seconds),
		xlabel=\papilo (time in seconds),
		xmode=log,
		ymode=log,
		grid=both,
		grid style={line width=.1pt, draw=gray!10},
		major grid style={line width=.2pt,draw=gray!50},
		minor tick num=5,
		xtick={0.01,0.1,1,10,100,1000},
		ytick={0.01,0.1,1,10,100,1000,10000},
		]
		\draw[domain=-100:100, smooth, variable=\x, black] plot ({\x}, {\x});
		\draw[domain=-100:100, smooth, variable=\x, black] plot ({\x}, {\x + 10});
		\addplot[scatter,only marks, scatter src=explicit symbolic]%
		table[meta=label] {
			x y label
			1000 1000 peer
			0.146 34.11 green
			0.046 10.56 green
			0.091 47.99 green
			3.98 2.15 green
			7.495 8372.310000000001 green
			0.616 2094.58 green
			1.147 3.86 green
			5.324 1942.23 green
			0.027 27.380000000000003 green
			0.203 116.2 green
			10.054 153.29 green
			0.018 6.4399999999999995 green
			1.192 47.099999999999994 green
			0.137 22.29 green
			0.037 14.83 green
			0.113 5.48 green
			2.705 8.49 green
			1.218 1966.6000000000001 green
			1.816 193.01 green
			0.027 51.93 green
			0.138 81.47 green
			0.036 7.06 green
			0.072 1.28 green
			0.554 37.199999999999996 green
			0.052 13.32 green
			0.018 2.27 green
			0.279 6.51 green
			0.623 140.26999999999998 green
			0.066 3.45 green
			1.696 3.3000000000000003 green
			11.227 2104.3300000000004 green
			0.256 117.19 green
			1.028 726.36 green
			0.009 3.3800000000000003 green
			9.253 2358.98 green
			0.38 42.29 green
			0.034 1.9899999999999998 green
			0.229 36.07 green
			0.507 1.06 green
			0.705 22.540000000000003 green
			0.444 29.84 green
			1.167 4087.21 green
			0.255 194.31 green
			0.432 133.09 green
			4.257 1.35 green
			0.067 4.89 green
			0.389 408.48999999999995 green
			4.288 3533.3500000000004 green
			0.271 578.04 green
			1.113 106.58 green
			0.365 14.17 green
			0.062 1.28 green
			0.004 3.99 green
			0.017 37.48 green
			0.153 2.0 green
			0.007 23.2 blue
			0.048 46.46 blue
			0.349 253.9 blue
			0.068 4.84 blue
			1.519 1.3 blue
			0.126 71.61 blue
			0.021 2.5 blue
			0.302 99.78999999999999 blue
			0.003 2.18 blue
			1.485 2104.29 blue
			0.309 2548.8199999999997 blue
			0.168 2.0100000000000002 blue
			2.851 2945.7000000000003 blue
			0.009 3.56 blue
			0.017 37.699999999999996 blue
			0.022 2.85 blue
			0.309 10.88 blue
			0.004 3.52 blue
			0.028 1.3599999999999999 blue
			0.229 35.81 blue
			2.845 2970.08 blue
			0.584 119.0 blue
			0.095 52.01 blue
			1.492 2097.74 blue
			0.096 52.48 blue
			0.049 46.66 blue
			2.498 8.48 blue
			0.134 33.62 blue
			5.02 4385.97 blue
			0.018 6.91 blue
			0.042 28.34 blue
			0.023 3.66 blue
			10.217 150.37 blue
			0.511 162.20000000000002 blue
			0.163 112.34 blue
			0.676 124.38 blue
			0.709 22.11 blue
			0.286 532.91 blue
			0.006 3.8200000000000003 blue
			0.299 10.92 blue
			0.046 181.37 blue
			0.038 15.79 blue
			0.006 3.8200000000000003 blue
			1.108 105.86 blue
			0.059 3.4000000000000004 blue
			0.509 160.9 blue
			0.059 3.91 blue
			1.086 105.59 blue
			0.315 2552.48 blue
			0.06 3.39 blue
			4.02 2.19 blue
			0.009 3.6 blue
			0.382 42.25 blue
			5.325 8491.69 blue
			0.143 112.64 blue
			0.134 71.98 blue
			0.149 2.0 blue
			0.095 2.15 blue
			0.051 1.27 blue
			1.541 1.31 blue
			66.524 606.9300000000001 blue
			0.038 15.809999999999999 blue
			1.095 3.87 blue
			0.801 44.21 blue
			1.11 3.86 blue
			5.378 8462.77 blue
			0.029 1.98 blue
			0.043 4.9799999999999995 blue
			0.01 1.22 blue
			0.678 124.2 blue
			0.115 5.3500000000000005 blue
			0.352 255.34 blue
			0.007 23.36 blue
			0.288 530.83 blue
			0.589 119.86 blue
			0.028 1.9899999999999998 blue
			0.199 125.76 blue
			0.042 28.23 blue
			0.231 36.07 blue
			0.897 886.9599999999999 blue
			0.152 1.9900000000000002 blue
			0.092 10.13 blue
			0.163 78.86 blue
			0.021 2.86 blue
			0.895 877.3499999999999 blue
			0.296 97.86999999999999 blue
			0.017 1.49 blue
			0.125 71.16 blue
			0.003 2.18 blue
			0.709 22.15 blue
			0.029 8179.29 blue
			7.315 4356.2300000000005 blue
			2.499 8.41 blue
			0.067 4.84 blue
			0.014 3.5300000000000002 blue
			7.29 4351.83 blue
			0.114 5.33 blue
			0.011 1.09 blue
			0.198 125.94 blue
			4.989 4367.419999999999 blue
			0.018 2.2600000000000002 blue
			0.018 6.87 blue
			0.017 2.08 blue
			0.845 39.19 blue
			0.542 1.06 blue
			0.004 3.4800000000000004 blue
			0.047 181.42 blue
			0.03 8214.94 blue
			
		};
		\end{axis}
		\end{tikzpicture}
	\end{minipage}
	\caption{
		Running times of \veripb vs. \papilo on test sets \pb (left) and \miplib (right), including all instances with more than 1~seconds in \veripb{} and less than 30~minutes in \papilo, and excluding timeouts.
		Green $+$ signs mark optimization and blue $\times$ signs mark decision instances.
	}
	\label{fig::experiment::PaPILOVeriPB}
\end{figure}

\subsection{Performance Analysis on Constraint Propagation}
\label{subsec::analysis-on-propagation}

Finally, we investigate how the performance of \veripb depends on whether we use RUP (as
in Sec.~\ref{subsec:impact-of-proof-logging-on-papilo} and Sec.~\ref{subsec::verifying-presolving}) or explicit cutting planes derivations (POL) to certify bound strengthening reductions from constraint propagation.
Here, we additionally exclude 9 large-scale instances\footnote{%
\solver{neos-4754521-awa\-rau.0.s}, \solver{neos-827015.0.s/u},
\solver{neos-829552.0.s/u}, \solver{s100.0.s/u},
\solver{normalized-datt256}, \solver{s100}}
for which \papilo reaches the memory limit when certifying with POL.
The results are summarized in Tab.~\ref{tab:rup_vs_add}.
The ``verified'' column contains the number of instances verified by \veripb within the time limit. The ``time'' column reports the time for verification.

Deriving the propagation directly with cutting planes 
is 3.2\% faster on \pb-dec, 2.8\% faster on \miplib-dec, 13.1\% faster on \miplib-opt, and 0.7\% faster on \pb-opt.
On 95\% of the decision instances using RUP is at most 9.7\% slower. 
While it is expected that verification is faster when the cutting planes proof is given explicitly, it is surprising that the performance difference between the methods is not more pronounced.
This is partly due to the cost of the watched-literal scheme~\cite{MMZZM01Engineering,SS06Pueblo} used by \veripb{} for unit propagation.
The overhead of maintaining the watches is present regardless of whether (reverse) unit propagation is used or not.
Furthermore, unit propagation is also used for automatically verifying redundance-based strengthening. 
Together, this limits the potential for runtime savings by providing the explicit cutting planes proof.

Furthermore, providing an explicit cutting planes proof for propagation requires printing the constraint into the certificate.
Hence, the certificate size becomes dependent on the number of non-zeros in the constraints leading to propagations.
In contrast, the overhead of RUP is constant and much smaller.

All in all, these results suggest to prefer RUP when deriving constraint propagation since it barely impacts the performance of \veripb{} and keeps the size of the certificate smaller.
The computational cost of RUP could be further reduced by extending it to accept an ordered list of constraints that shall be propagated first, similar as in~\cite{LRAT}.
Such an extension could also be used for other presolving techniques, in particular probing and simple probing.


\begin{table}[t]
	\centering
	\small
	\caption{Comparison of the runtime of \veripb with RUP and POL over instances with at least 10 propagations.}
	\begin{tabular*}{\textwidth}{@{}l@{\;\;\extracolsep{\fill}}lrrrrr}
		\toprule
		& &  \multicolumn{2}{c}{RUP} & \multicolumn{2}{c}{POL} &  \\[-0.5ex] 
		\cmidrule{3-4}
		\cmidrule{5-6}
		test set & size &  verified & time [s] & verified & time [s] & relative \\[-0.5ex]
		\midrule
		\pb-dec & 284 & 284 & 2.21 & 284 & 2.14 &  0.968\\[-0.5ex]
		\miplib-dec & 35 & 31 & 153.23 & 31 & 148.88 &  0.972\\[-0.5ex]
		\pb-opt & 153 & 142 & 28.43 & 142 & 28.22 &  0.993\\[-0.5ex]
		\miplib-opt & 16 & 14 & 147.11 & 14 & 127.83 &  0.869\\[-0.5ex]

		\bottomrule
	\end{tabular*}
	\label{tab:rup_vs_add}
\end{table}

\section{Conclusion}
\label{sec::conclusion}

In this paper we set out to demonstrate how presolve techniques from state-of-the-art MIP solvers can be equipped with certificates in order to verify the equivalence between original and reduced models.
Although the pseudo-Boolean proof logging format behind \veripb{}~\cite{BGMN22Dominance} was not designed with this purpose in mind, we could show that a limited extension needed for handling updates of the objective function is sufficient to craft a certified presolver for \mbox{$0$--$1$} ILPs.

However, our experimental study on instances from pseudo-Boolean competitions and \miplib also exhibited that the verification of MIP-based presolving can suffer from large and overly verbose certificates.
To shrink the proof size we introduced a sparse objective update function but identified further possible improvements.
%
First, a native substitution rule in \veripb{} would remove the need for the explicit
  derivation of new aggregations and the verification of checked deletion as described in
  Sec.~\ref{subsec::general}. For instances where presolving is dominated by substitutions, we
  estimate that this would reduce certificate sizes by up to 90\%, and no more
  time would be spent on checked deletion for substitutions.
Second, augmenting the RUP syntax by the option to specify an ordered list of constraints to propagate first, similarly as in~\cite{LRAT}, would accelerate RUP, in particular for fast verification of bound strengthenings by constraint propagation.

While \veripb{} is currently restricted to operate on integer coefficients only, the certification techniques presented in Sec.~\ref{sec::certifying} do not rely on this assumption and are applicable to general binary programs.
It has been shown how to construct \veripb{} certificates for bounded integer domains~\cite{GMN22AuditableCP,MM23ProofLogging}, and within the framework of the generalized proof system laid out in~\cite{DEGH23ProofSystem}, our certificates would even translate to continuous and unbounded integer domains.
To conclude, we believe our results show convincingly that this type of proof logging
techniques is a very promising direction of research also for MIP presolve beyond \zeroone
ILPs.

\clearpage
\medskip\noindent{\small
\textbf{Acknowledgements.}
The authors wish to acknowledge helpful technical discussions on
\veripb in general and the objective update rule in particular with
Bart Bogaerts, Ciaran McCreesh, and Yong Kiam Tan.
The work for this article has been partly conducted within the Research Campus MODAL funded by the German Federal Ministry  of Education and Research (BMBF grant number 05M14ZAM). 
Jakob Nordstr\"om was supported by the Swedish Research Council grant 2016-00782 and the Independent Research Fund Denmark grant 9040-00389B.
Andy Oertel was supported by the Wallenberg AI, Autonomous Systems and Software Program (WASP) funded by the Knut and Alice Wallenberg Foundation.
The computational experiments were enabled by resources provided by LUNARC at Lund University.

}

\bibliographystyle{splncs04}
\bibliography{bibliography,refArticles,refBooks,refOther}

\end{document}